\documentclass[12pt]{article}
\usepackage[T1]{fontenc} 
\usepackage{hyperref}
\usepackage{amsmath}
\usepackage{amsfonts}
\usepackage{amssymb}

\begin{document}

\title{Tr\`es courte enqu\^ete sur l'extension non-triviale de la logique de propositions \`a la logique du premier et deuxi\`eme ordre\footnote{Short essay written in 2005 in the context of a Master's course in philosophy and classical logic.}}

\author{Hector Zenil\footnote{\href{mailto:hector.zenil-chavez@malix.univ-paris1.fr}{hector.zenil-chavez@malix.univ-paris1.fr}}\\
Institut d'Histoire et de Philosophie des Sciences et des Techniques\\
(Paris 1 Panth\'eon-Sorbonne/ENS Ulm/UMR 8590 CNRS)}

\date{}

\maketitle

\begin{abstract}
La construction formelle de la logique de second ordre ou calcul de pr\'edicats consiste seulement \`a ajouter essentiellement de quantificateurs \`a la logique des propositions. Pourquoi la logique de deuxi\`eme ordre ne peut pas se r\'eduire \`a celle de premier ordre ? Comment d\'emontrer que certains pr\'edicats sont d'ordre sup\'erieur ? Quel type d'ordre corresponde au langage naturel ? Est-ce qu'il y a une position philosophique derri\`ere chaque logique, m\^{e}me pour les logiques classiques ? Quelle position pour quelle puissance ? Ce sont des questions qu'on se pose et qu'on esquisse.
\end{abstract}

\section{Contexte}

\noindent La construction formelle de la logique de premi\`ere ordre ou calcul de pr\'edicats consiste \`a ajouter essentiellement de quantificateurs \`a la logique des propositions. Ces nouveaux quantificateurs (qui l'on peut \^{e}tres repr\'esent\'es par un seul car on arrive \`a l'autre par moyen de la n\'egation) quantifient un deuxi\`eme type d'objets qui peuvent \^{e}tre des pr\'edicats ou des relations. Ceux relations unaires peuvent se voir aussi comme des ensembles. D'un point de vue s\'emantique, cela revient \`a dire que l'on consid\`ere les fonctions et les pr\'edicats comme des objets \`a part enti\`ere, au m\^{e}me titre que, par exemple, les variables des nombres.

\section{Pourquoi la logique de deuxi\`eme ordre ne peut pas se r\'eduire \`a celle de premier ordre}

On peut penser qu'il suffit d'ajouter au domaine d'un langage de premier ordre les \'el\'ements sur lesquels les quantificateurs de la logique de deuxi\`eme ordre quantifient. Cependant, m\^{e}me si c'est possible, le syst\`eme r\'esultant est moins puissant que l'original puisque tous les ensembles ne sont pas d\'efinis par une formule. On sait que la logique de deuxi\`eme ordre a des propri\'et\'es diff\'erentes qui font une diff\'erence essentielle entre elles. Il s'agit des propri\'et\'es de compl\'etude: la descriptive, la s\'emantique et la syntaxique. 

\subsection{Compl\'etude descriptive}

Une grande partie des math\'ematiques peut \^{e}tre formalis\'ee en logique du premier ordre gr\^{a}ce \`a l'\'egalit\'e. Par exemple, la th\'eorie des groupes. Cependant, il existe des notions qui ne peuvent pas \^{e}tre saisies pr\'ecis\'ement par aucune th\'eorie du premier ordre comme la propriet\'e math\'ematique d'\^{e}tre infini ou d'\^{e}tre un nombre naturel. Nous disons que de telles notions ne sont pas axiomatisables au premier ordre. Par exemple il n'y a pas d'axiomatisation au premier ordre dont le seul mod\`ele soit l'ensemble $N$ des entiers naturels muni des op\'erations arithm\'etiques usuelles.

On peut construire une arithm\'etique, l'arithm\'etique de Peano ($PA$), mais il faudra un nombre infini d'axiomes, autrement elle le calcul de pr\'edicats ne d\'ecrira pas compl\`etement l'arithm\'etique car le sch\'ema de r\'ecurrence est une collection infinie d'axiomes. Si on \'evite l'axiome de r\'ecurrence on produit une aritm\'etique moins puissante, l'arithm\'etique de Robinson. En fait le sch\'ema de r\'ecurrence a besoin d'une logique plus puissante, une logique de deuxi\`eme ordre. On ne peut pas saisir le sch\'ema dans la logique du premier ordre parce qu'on ne peut pas quantifier que sur des variables, et les variables d\'enotent des valeurs, pas des pr\'edicats sur les valeurs. On remplace couramment donc l'axiome de r\'ecurrence par une version plus faible. L'effet secondaire de ne pas pouvoir saisir compl\'etement $PA$ est que'elle \'a d'autres mod\`eles, les mod\`eles non-standard de l'arithm\'etique dans le calcul de pr\'edicats. 

Le domaine de la compl\'etude descriptive s'agit donc de pouvoir expressif pour saisir de mani\`ere finie, par exemple $PA$. La logique de deuxi\`eme ordre r\'ealise cet objective et donc elle est plus puissante dans ce sens que la logique de premier ordre, car son pouvoir expressif est plus grand.

\subsection{Compl\'etude s\'emantique}

Il s'agit de l'impossibilit\'e d'avoir un syst\`eme de d\'eduction comme celui de la logique de premier ordre. Il y a des mod\'eles de deuxi\`eme ordre valides qui ne sont pas d\'eriv\'es par le syst\`eme de d\'eduction. Dans les syst\`emes de premier ordre, on peut donner un nombre fini d'axiomes, de sch\'emas d'axiomes ou de r\`egles de d\'eduction tels que toutes les preuves logiquement valides formul\'ees avec la grammaire du calcul des pr\'edicats du premier ordre soient obtenues \`a partir de ces principes. Par exemple, le syst\`eme de d\'eduction \`a la Hilbert o\`u les quinze r\`egles de d\'eduction naturelle sont un syst\`eme compl\`et des principes de la logique du premier ordre. Cela veut dire que si $p$ et $q$ sont deux formules du calcul des pr\'edicats, on dira que $q$ est une cons\'equence logique de $p$ lorsque tout mod\`ele de p est aussi un mod\`ele de $q$, autrement dit, pour tout mod\`ele m, si $p$ est vrai dans m alors $q$ est vrai dans $m$. 

$D$ est incomplet, autrement dit s'il existe une formule $F$ construite au moyen des connecteurs logiques, et des quantificateurs sur les formules tel que $F$ soit vrai dans le mod\`ele standard $N$ de l'arithm\'etique, mais n'est pas prouvable dans $D$.

Cependant, l'ensemble de toutes les formules universellement valides de deuxi\`eme ordre n'est pas \'enum\'erable m\^{e}me r\'ecursivement. Ceci est une cons\'equence du th\'eor\`eme de l'incompl\'etude. Et donc quel que soit les tentatives pour d\'efinir un syst\`eme de d\'eduction, ce ne sera pas suffisant pour qu'il soit complet par rapport aux formules valides.

On sait depuis le travail de G\"odel qu'une th\'eorie formelle contenant l'arithm\'etique n'est pas compl\`ete et qu'elle contient des \'enonc\'es vrais ind\'emontrables.

\subsection{Compl\'etude syntaxique}

Un syst\`eme formel est dit complet syntactiquement si, en utilisant ses axiomes et ses r\`egles, on peut montrer soit qu'il est un th\'eor\`eme du syst\`eme ou que la n\'egation est un th\'eor\`eme.

\section{Les propri\'et\'es de la logique de deuxi\`eme ordre}

On a vu qu'avec la logique de deuxi\`eme ordre, on gagne en compl\'etude descriptive et en m\^{e}me temps qu'on perd les propri\'et\'es les plus importantes de la logique de premier ordre. 

D'ailleurs, \'etant donn\'e la taille du domaine des syst\`emes de deuxi\`eme ordre, la complexit\'e de calcul est beaucoup plus importante. Elle est intraitable et sa validit\'e ne peut pas \^{e}tre toujours calculable par une machine de Turing (seulement avec une machine oracle car les sous-ensembles infinis du domaine d'un syst\`eme de deuxi\`eme ordre peuvent se voir comme les expansions d\'ecimales des nombres r\'eels et nous savons que la plupart des nombres r\'eels ont un degr\'e de Turing mais il ne sont pas toujours dans la hi\'erarchie arithm\'etique, c'est \`a dire qu'ils n'ont pas d'expression sous forme de formule).

Pour v\'erifier qu'on gagne en pouvoir expressif, on peut voir qu'il suffit de la logique de deuxi\`eme ordre pour cat\'egoriser le mod\`ele standard de l'arithm\'etique et que tous les mod\`eles non-standards de l'arithm\'etique de deuxi\`eme ordre sont isomorphes et donc il y a seulement une arithm\'etique de deuxi\`eme ordre contrairement \`a ce qui se passe dans le calcul de pr\'edicats qui admet trop de mod\`eles non isomorphes.

\section{Comment d\'emontrer que certains pr\'edicats sont de deuxi\`eme ordre}

Dire qu'un pr\'edicat est de deuxi\`eme ordre signifie qu'il ne peut pas \^{e}tre \'ecrit comme un \'enonc\'e de premier ordre. On conna\^{i}t quelques exemples fournis par Quine. Pour d\'emontrer qu'un \'enonc\'e ne peut pas \^{e}tre \'ecrit dans le premier ordre, il suffit de voir qu'il n'est pas vrai dans l'interpr\'etation standard de l'arithm\'etique de premier ordre mais qu'il est vrai dans l'arithm\'etique d'ordre sup\'erieur. Pour r\'eussir cela, il suffit d'arriver \`a une contradiction banale concernant les successeurs puisque tout nombre non-standard est un successeur de quelque autre tandis que dans le mod\`ele non-standard le z\'ero n'est pas le successeur d'aucun autre nombre. Voici un exemple:
\\

\noindent (1): ``Soient certains pistoleros. Chacun d'eux a tir\'e sur le pied droit d'au moins l'un d'entre eux" (There are some gunslingers each of whom has shot the right foot of at least one of the others).
\\

\noindent Qui peut \^{e}tre \'ecrit en calcul de pr\'edicats comme:
\\

\noindent (2): $\exists$$X$($\exists$$x$$X$$x$$\wedge$$\forall$$x$($X$$x$$\rightarrow$$\exists$$y$($X$$y$$\wedge$$y$$\neq$$x$$\wedge$$B$$x$$y$)))
\\

\noindent Il s'agit alors de trouver une substitution de telle sorte que l'\'enonc\'e soit vrai dans le mod\`ele non-standard et faux dans les mod\`eles standards.
\\

\noindent On fait la substitution $x$$=$$y$$+$$1$ dans la relation $B$$x$$y$, car il s'agit toujours de se profiter de la caract\'eristique des nombres non-standards, qui ont la particularit\'e que si a est un nombre non-standard, alors son pr\'ed\'ecesseur est aussi un nombre non-standard. C'est \`a dire si $a$ est non-standard, et $a$$=$$b$$+$$1$, alors $b$ est non-standard. Donc le z\'ero n'est pas un nombre non-standard.\\

\noindent Donc l'\'enonc\'e (2) devient :
\\

\noindent (3): $\exists$$X$($\exists$$x$$X$$x$$\wedge$$\forall$$x$($X$$x$$\rightarrow$$\exists$$y$($X$$y$$\wedge$$y$$\neq$$x$$\wedge$$x$$=$$y$$+$$1$)))
\\

\noindent On d\'emontre que :
\\

\noindent $M$$\models$(2) pour $M$ non-standard n'est pas vrai dans le cas de $M'$$\models$(2) pour $M'$ standard.
\\

\noindent Supposons que $M'$$\models$(2) pour $M'$ standard quand  l'ant\'ec\'edent \\ ``$\exists$$X$($\exists$$x$$X$$x$$\wedge$$\forall$$x$($X$$x$))'' est vrai. Alors le cons\'equent ``$\exists$$y$($X$$y$$\wedge$$y$$\neq$$x$$\wedge$$x$$=$$y$$+$$1$)'' devrait \^{e}tre vrai. 
\\

\noindent Cependant si $\exists$$y$($X$$y$$\wedge$$y$$\neq$$x$) et $x$$=$$y$$+$$1$ $y$$=$$x$$-$$1$. C'est \`a dire que $y$ a un pr\'ed\'ecesseur, ce qui est faux dans le cas $y$$=$$0$. Donc contradiction.
\\

\noindent Donc $M'$$\models$(2) pour $M'$ standard est faux.
\\

\noindent Tandis que pour les nombres non-standards, il est vrai. Donc $M$$\models$(2) pour $M$ non-standard.
\\

\noindent Alors (3) est effectivement un \'enonc\'e de deuxi\`eme ordre, donc (2) et (1) aussi.

\section{Une position philosophique derri\`ere chaque logique, m\^{e}me pour les logiques classiques ?}

Il semble que le calcul de pr\'edicats comme syst\`eme de expression math\'ematique ne soit pas $philosophiquement$ $neutre$ car elle peut donner lieu \`a certaines interpr\'etations non r\'ealistes et \`a de mod\`eles dites non-standards (qu'on n'avait pas l'intention de produire initialement). Quel sont les compromis entre expressivit\'e, puissance et r\'ealisme ? On pourrait dire si la logique de second ordre est plus ou moins r\'ealiste ? et dans quel sens ?

\end{document}